\documentclass[11pt,twoside]{article}
\usepackage{latexsym,amsmath}
\usepackage{indentfirst}
\usepackage{graphicx}

\numberwithin{equation}{section}
\newtheorem{Prop}{\bf Proposition}[section]

\newtheorem{Rem}{\bf Remark}[section]
\newtheorem{Ex}{\bf Example}[section]

\begin{document}
\def \b{\Box}

\begin{center}
{\Large {\bf Dynamics analysis of the fractional-order Lagrange system}}
\end{center}

\begin{center}
{\bf Mihai IVAN}
\end{center}

\setcounter{page}{1}

\pagestyle{myheadings}

{\small {\bf Abstract}. The main purpose of  this paper is to study the fractional-order model with Caputo derivative associated to Lagrange system. For this fractional-order system we investigate the existence and uniqueness of solutions of initial value problem, asymptotic stability of its equilibrium states,  stabilization problem using appropriate controls and numerical integration via the fractional Euler method.}\\[-0.4cm]
{}\\[-0.5cm]

\section {Introduction}
\smallskip
\indent

The dynamic behavior of differential systems (in particular, Hamilton-Poisson systems) have been studied  due to their deep applications in many
areas of science and engineering including robotics, spatial dynamics, complete synchronization and secure comunications.
 Three remarkable classes of systems formed by differential equations on ${\bf R}^{3}$ are the general Euler top system, the  systems of Maxwell-Bloch type and the family of Lotka-Volterra systems. In the recent decades, a series of dynamical systems belonging to these three classes have been studied from the point of view of Poisson geometry by many researchers, for instance  \cite {gunu, takh, puta, pagi, igmp, giv1, lima,opmi, enpp, laph}.

An interesting example of a differential system of the Euler top type is the {\it Lagrange system} \cite{takh}. It is described by the following differential equations on $~{\bf R}^{3}:$\\[-0.2cm]
\begin{equation}
 \dot{x}^{1}(t)  =
 x^{2}(t) x^{3}(t),~~~
 \dot{x}^{2}(t)  =
  x^{1}(t)x^{3}(t),~~~
  \dot{x}^{3}(t) = x^{1}(t)x^{2}(t),\label{(2.1)}
\end{equation}
 where $ \dot{x}^{i}(t) = d x^{i}(t)/dt $ and  $ t $ is the time.

The system  $ (1.1) $ is used in theoretical physics for the study of $ SU(2)-$ monopoles.

The fractional calculus deals with derivations and integration of arbitrary order and has deep and natural connections with many fields of applied mathematics, mathematical physics, engineering, biological systems, control processing, chaos synchronization, \cite {mati, podl, agra, ahma, nabu, odmo, diet, bhad, badi, danc, ahme}.
The fractional models associated to nonlinear dynamical systems have been investigated in the papers \cite {migo, giv2, ivmi}.

The Lie groups, Lie algebroids and Leibniz algebroids have proven to be powerful tools for geometric formulation of the Hamiltonian mechanics \cite{demi,giop, gimo}. Also, they have been used  in the investigation of many fractional dynamical systems \cite{nabu, imod, migo, migi}.

This paper is structured as follows. The fractional-order Lagrange system is defined in Section 2. The existence and uniqueness of solutions of initial value problem for the  fractional model $(2.2)~$ is discussed.  For the asymptotic stabilization problem of fractional model  $~(2.2), $ we associate the fractional Lagrange system with  controls $~c_{1} $ and $ c_{2},~$ denoted by $ (3.3).~$ In Propositions $~(3.1)-(3.3)~$ are established sufficient conditions on parameters  $~c_{1} $ and $ c_{2}~$  to control the chaos in the fractional-order Lagrange system $~(3.3).~$  In Section 4, we give the numerical integration of fractional-order Lagrange system $(3.3),$ via the fractional Euler's method. Finally, a numerical simulation of fractional model  $~(3.3)~$ is shown to verify the theoretical results.\\[-0.5cm]

\section{ The fractional-order Lagrange system}

We recall the Caputo definition of fractional derivatives, which
is often used in concrete applications. Let $ f\in C^{\infty}(
\textbf{R}) $ and $ q \in \textbf{R}, q
> 0. $ The $ q-$order Caputo differential operator
\cite{diet}, is described by $~D_{t}^{q}f(t) = I^{m -
q}f^{(m)}(t), ~q > 0,~$ where $~f^{(m)}(t)$  represents the $
m-$order derivative of the function $ f,~m \in \textbf {N}^{\ast}$
is an integer such that $ m-1 \leq q \leq m $ and $ I^{q} $ is the
{\it $ q-$order Riemann-Liouville integral operator} \cite{podl, diet}, which
is expressed by:\\[-0.4cm]
\begin{equation}
I_{t}^{q}f(t) =\displaystyle\frac{1}{\Gamma(q)}\int_{0}^{t}{(t-s)^{q
-1}}f(s,y(s))ds,~q > 0,~~~s\in [0, T],  \label{(2.1)}
\end{equation}
where $~\Gamma $ is the Euler Gamma  function. If $ q =1$, then $ D_{t}^{q}f(t) = df/dt.$
\markboth{M. Ivan}{Dynamics analysis of the fractional-order Lagrange system}

The {\it fractional-order Lagrange system}
associated to dynamics $ (1.1)$ is defined by the following
set of fractional differential equations:\\[-0.2cm]
\begin{equation}
\left\{ \begin{array} {lcl}
 D_{t}^{q}{x}^{1}(t) & = &
 x^{2}(t) x^{3}(t)  \\[0.1cm]
 D_{t}^{q}{x}^{2}(t) & = & x^{1}(t) x^{3}(t),~~~ q \in (0,1),\\[0.1cm]
  D_{t}^{q}{x}^{3}(t) & = & x^{1}(t)x^{2}(t). \label{(2.1)}
  \end{array}\right.
\end{equation}
The initial value problem of the fractional-order Lagrange system
$(2.1)$ can be represented in the following matrix
form:\\[-0.2cm]
\begin{equation}
D_{t}^{\alpha}x(t)  =  x^{1}(t) A x(t) + x^{3}(t) B  x(t) ,~~~~~ x(0) =
x_{0},\label{(2.2)}
\end{equation}
where $0 < q < 1,~ x(t)= ( x^{1}(t),
 x^{2}(t), x^{3}(t))^{T}, ~t\in(0,\tau)$ and\\[-0.1cm]
\[
A = \left ( \begin{array}{ccc}
0 & 0 & 0 \\
0 & 0 & 0 \\
0 & 1 & 0 \\
\end{array}\right ),~~~ B = \left ( \begin{array}{ccc}
0 & 1 & 0 \\
1 & 0 & 0\\
0 & 0 & 0\\
\end{array}\right ).
\]
\begin{Prop}
The initial value problem of the fractional-order Lagrange system
$(2.1)$ has a unique solution.
\end{Prop}
{\it Proof.} Let $ f(x(t))= x^{1}(t) A x(t) + x^{3}(t) B x(t).$ It is
obviously continuous and bounded on $ D =\{ x \in {\bf R}^{3} |~
x^{i}\in [x_{0}^{i} - \delta, x_{0}^{i} + \delta]\}, i=\overline{1,3} $ for any
$\delta>0. $ We have $~f(x(t)) - f(y(t)) =  x^{1}(t) A x(t) - y^{1}(t) A y(t) + x^{3}(t) B x(t)-  y^{3}(t) B y(t)=g(t)+h(t), $
where $~g(t)= x^{1}(t) A x(t) - y^{1}(t) A y(t)~$ and $~h(t)= x^{3}(t) B x(t) - y^{3}(t) B y(t). $
 Then\\[0.2cm]
$(a)~~|f(x(t)) - f(y(t))|\leq |g(t)| + |h(t)|. $

It is easy to see that $~g(t)= (x^{1}(t)- y^{1}(t))A
x(t)+ y^{1}(t) A(x(t)- y(t)).~$ Then\\[0.2cm]
$|g(t)| \leq |(x^{1}(t)- y^{1}(t)) A x(t)| + |y^{1}(t)
A (x(t)- y(t))|\leq $\\
$\leq  \|A\|( |x(t)|\cdot |x^{1}(t)- y^{1}(t)| + |y^{1}(t)|\cdot |x(t) -y(t)|, $\\
where $ \|\cdot \| $ and $|\cdot|$ denote matrix norm
and vector norm respectively.

Using the inequality $~|x^{1}(t)- y^{1}(t))|\leq |x(t)- y(t))|~$ one obtains\\[0.1cm]
$(b)~~~|g(t)| \leq (\|A\|+ |y^{1}(t)|)\cdot |x(t)- y(t)|.$\\[0.1cm]
Similarly, we have\\[0.2cm]
$(c)~~~|h(t)| \leq (\|B\| + |y^{3}(t)|)\cdot |x(t)- y(t)|.$\\[0.1cm]
According to $(b)$ and $ (c),$ the relation $(a)$ becomes\\[0.1cm]
$(d)~~~|f(x(t)) - f(y(t))|\leq  (\|A\| + \|B\| +
|y^{1}(t)| + |y^{3}(t)|)\cdot |x(t)-y(t)|.$\\[0.1cm]
Replacing  $\|A\|= 1,~ \|B\| = \sqrt{2}$ and using the inequalities $~|y^{i}(t)|\leq |x_{0}| + \delta, ~i=1,3~$ from the
relation $ (d), $ we deduce that\\[0.1cm]
$(e)~~~|f(x(t)) - f(y(t))|\leq  L\cdot |x(t)-y(t)|,~~~~~
\hbox{where}~ L =1 + \sqrt{2} + 2(|x_{0}| + \delta) > 0.$

 The inequality $(e)$ shows that $ f(x(t))$
satisfies a Lipschitz condition. Based on the results of Theorems
$1$ and $2$ in \cite{difo}, we can conclude that the initial value
problem of the system $(2.2)$ has a unique solution. \hfill$\Box$

For the fractional Lagrange system $(2.1)$ we introduce the following
 notations:\\[-0.2cm]
 \begin{equation}
f_{1}(x) = x^{2}  x^{3},~~~ f_{2}(x) =
x^{1} x^{3},~~~ f_{3}(x)  = x^{1} x^{2}.\label{(2.2)}
\end{equation}
\begin{Prop}
{\it The equilibrium states of the fractional Lagrange system
$(2.1)$ are given as the union of the following three
families:\\[-0.3cm]
\[
 E_{1}:=\{ e_{1}^{m}=(m, 0,
0)\in {\bf R}^{3} |~m \in {\bf R}
\},~~~ E_{2}:=\{ e_{2}^{m}=(0, m, 0)\in {\bf R}^{3} |~ m \in {\bf
R}\},
\]
$ E_{3}:=\{ e_{3}^{m}=(0, 0, m) \in {\bf R}^{3}  |~ m \in {\bf
R}\}.$}
\end{Prop}
{\it Proof.} The equilibrium states are solutions of the equations
$~f_{i}(x)=0, i=\overline{1,3}$ where $~f_{i},~i=\overline{1,3}$
are given by (2.2).\hfill$\Box$

Let us we present the study of asymptotic stability of equilibrium
states for the fractional system $(2.1)$. Finally, we will discuss
how to stabilize the unstable equilibrium states of the system
$(3.1)$ via fractional order derivative. For this study we apply
the Matignon's test.

The Jacobian matrix associated to system $(2.1)$ is:
\[
J(x) = \left ( \begin{array}{ccc}
0 & x^{3}    & x^{2} \\
  x^{3}   & 0 &  x^{1}\\
   x^{2}   &  x^{1}  & 0 \\
\end{array}\right ).\\[-0.1cm]
\]
\begin{Prop} {\rm (\cite{mati})}
Let $x_{e}$ be an equilibrium state of system $(2.1)$ and
$ J(x_{e})$ be the Jacobian matrix $J(x)$ evaluated at
$x_{e}. $ Then  $~ x_{e}$ is locally asymptotically stable, iff all eigenvalues of the matrix $J(x_{e})$
 satisfy the inequality: $~~| arg(\lambda (J_(x_{e}))) | > \displaystyle\frac{q \pi}{2}.$ \hfill$\Box$\\[-0.2cm]
\end{Prop}
\begin{Prop}
The equilibrium states $ e_{i}^{m}\in E_{i},~i=\overline{1,3}~$
are unstable $ (\forall) q \in (0,1).$
\end{Prop}
{\it Proof.} The characteristic polynomial of the matrix $~
J(e_{1}^{m}) =\left (\begin{array}{ccc}
  0 & 0  & 0\\[0.1cm]
  0 & 0 & m \\
  0 & m  & 0 \\
\end{array}\right ) $
is\\
 $~ p_{J(e_{1}^{m})}(\lambda) = \det ( J(e_{1}^{m}) -
\lambda I) = - \lambda ( \lambda ^{2} -  m^{2}).~$  For $ m = 0, $ the characteristic polynomials of  $ J(e_{0}) $ is
  $~p_{J(e_{0})}(\lambda) = \lambda ^{3}.$

The characteristic polynomials of  $
J(e_{i}^{m}),~  i=2,3 $ are $ p_{J(e_{i}^{m})}(\lambda) = - \lambda ( \lambda^{2} - m^{2} ).$

 The equations $~p_{J(e_{0})}(\lambda) = 0 $ and  $ ~p_{J(e_{i}^{m})}(\lambda) = 0 $ for $ i= \overline {1,3} $ have
 the root $ \lambda_{1} = 0 .$  Since $~arg(\lambda_{1}) =0 < \frac{q \pi}{2}$
for all $ q\in (0,1),$  by Proposition 2.3 follows that the equilibrium states $ e_{0} $ and $ e_{i}^{m}, i=\overline{1,3}$ are unstable for
all $ q\in (0,1).$ \hfill$\Box$\\[-0.5cm]

\section{ Stability analysis of the fractional-order Lagrange system with controls $~(3.3)~$}

In the case when $x_{e}$ is a unstable equilibrium state of the
fractional-order system $(2.1)$, we associate to $(2.1)$ a new
fractional system, called the {\it fractional-order Lagrange
system with controls} and given by:\\[-0.2cm]
\begin{equation}
\left\{ \begin{array} {lcl}
 D_{t}^{q}{x}^{1}(t) & = &
 x^{2}(t) x^{3}(t) + u_{1}(t) \\[0.1cm]
 D_{t}^{q}{x}^{2}(t) & = & x^{1}(t) x^{3}(t) + u_{2}(t) ,~~~~~~~~~~~~~~~~~ q \in (0,1),\\[0.1cm]
  D_{t}^{q}{x}^{3}(t) & = &  x^{1}(t)x^{2}(t) + u_{3}(t), \label{(3.1)}
  \end{array}\right.
\end{equation}
where $ u_{i}(t), i=\overline{1,3}$ are control functions.\\

In this section we take the control functions $ u_{i}(t), i=\overline{1,3}, $ given by:\\[-0.2cm]
 \begin{equation}
u_{1}(t) = c_{1} x^{1}(t),~~~ u_{2}(t) = c_{2}
x^{2}(t),~~~ u_{3}(t)  = c_{2} x^{3}(t), ~~~ c_{1}, c_{2} \in {\bf R}.\label{(3.2)}
\end{equation}
With the control functions $(3.2), $ the system $(3.1) $ becomes:\\[-0.2cm]
\begin{equation}
\left\{ \begin{array} {lcl}
 D_{t}^{q}{x}^{1}(t) & = &
 x^{2}(t) x^{3}(t) + c_{1} x^{1}(t) \\[0.1cm]
 D_{t}^{q}{x}^{2}(t) & = & x^{1}(t) x^{3}(t) + c_{2}x^{2}(t) ,~~~~~~~~~~~~~~~~~ q \in (0,1),\\[0.1cm]
  D_{t}^{q}{x}^{3}(t) & = &  x^{1}(t) x^{2}(t) + c_{2} x^{3}(t) \label{(3.3)}
  \end{array}\right.
\end{equation}
where $ c_{1}, c_{2} \in {\bf R}^{\ast} $  are real constants.\\

The system $(3.3) $ is called the {\it fractional-order Lagrange system with controls  $ c_{1}, c_{2}. $}

If one selects the parameters $ c_{1}, c_{2} $ which then make the eigenvalues of the Jacobian matrix of fractional model
$(3.3)$ satisfy one of the conditions from Proposition 2.3, then its trajectories asymptotically approaches the unstable
equilibrium state $ x_{e} $ in the sense that $\lim_{t\rightarrow \infty} \|x(t)- x_{e}\|= 0$, where $\|\cdot\|$ is the Euclidean norm.

The Jacobian matrix of the fractional-order Lagrange system  $(3.3)$ with  controls $ c_{1}, c_{2} $ is:\\[-0.2cm]
\[
J(x, c_{1}, c_{2}) = \left ( \begin{array}{ccccc}
c_{1} &  x^{3} & x^{2} \\
  x^{3}  & c_{2} &   x^{1}\\
   x^{2}   &  x^{1}  & c_{2} \\
\end{array}\right ).
\]
\begin{Prop}
Let  be the fractional-order Lagrange system $(3.3).$

$(i)~$ If  $~c_{1} < 0, c_{2} < 0, $ then $ e_{0}$ is asymptotically stable $~(\forall)~q\in (0,1);$

$(ii)~$ If  $~c_{1} > 0, c_{2} > 0 $ or $~c_{1} c_{2} < 0, $ then $ e_{0}$ is unstable $~(\forall)~q\in (0,1).$
\end{Prop}
{\it Proof.} The characteristic polynomial of $~J(e_{0}, c_{1}, c_{2}) $ is $~p_{0}(\lambda) = -(\lambda - c_{1}) (\lambda - c_{2}) ^{2}.$
The roots of the equation $ p_{0}(\lambda)= 0~$ are $\lambda_{1}=
c_{1},~\lambda_{2,3}= c_{2}.$\\
$(i)~$ We suppose  $~c_{1} < 0~$ and $~c_{2} < 0.$ In this case  we have $ Re (\lambda_{i})<
0 $ for $ i=\overline{1,3}. $ Since $|arg(\lambda_{i})| =\pi > \displaystyle\frac{q \pi}{2},
i=\overline{1,3} $ for all $ q\in (0, 1)$, by Proposition 2.3, it implies that $e_{0}^{m}$ is locally asymptotically stable.\\
$(ii)~$ We suppose  $~c_{1} > 0, c_{2} > 0 $ or $~c_{1} c_{2} < 0. $  Since
$ J(e_{0}, c_{1}, c_{2})$ has at least a positive eigenvalue, it follows that  $ e_{0} $ is unstable.  Hence,  $ (i) $ and $(ii) $ hold. \hfill$\Box$
\begin{Prop}
Let  be the fractional-order Lagrange system $(3.3)~$ and $~q\in (0,1).$\\
$(i)~$ If  $~c_{1} < 0, c_{2} < 0, $ then $~e_{1}^{m} $ is asymptotically stable $~(\forall) m\in ( c_{2}, - c_{2})~$ and unstable $~(\forall)~ m\in ( -\infty, c_{2}) \cup ( - c_{2}, \infty).$\\
$(ii)~$ If  $~c_{1} > 0, c_{2} > 0 $ or $~c_{1} c_{2} < 0, $ then $ e_{1}^{m} $ is unstable $~(\forall) m\in {\bf R}^{\ast}.$
\end{Prop}
{\it Proof.} The Jacobian matrix of  $ (3.3)$ at $ e_{1}^{m}$ is
$~J(e_{1}^{m}, c_{1}, c_{2})=\left (
\begin{array}{ccccc}
c_{1} & 0 &  0\\[0.2cm]
0 & c_{2} & m \\
0 & m  &  c_{2}\\
\end{array}\right ) $
whose characteristic polynomial is $~p_{1}(\lambda) = -(\lambda -
c_{1})[ (\lambda - c_{2}) ^{2} - m^{2}].$
The roots of the equation $ p_{1}(\lambda)= 0~$ are $\lambda_{1}=
c_{1},~\lambda_{2,3}= c_{2} \pm m.$\\
$(i)~$ We suppose  $~c_{1} > 0, c_{2} > 0. $ The eigenvalues $~\lambda_{i},~i=\overline{1,3}$ are all negative
if and only if  $~m\in ( c_{2}, - c_{2}). $
  Then $ e_{1}^{m}$ is asymptotically stable.\\
$(ii)~$ We suppose  $~c_{1} > 0, c_{2} > 0 $ or $~c_{1} c_{2} < 0. $ In these cases,the matrix
$ J(e_{0}, c_{1}, c_{2})$ has at least a positive eigenvalue. It follows that  $ e_{1}^{m} $ is unstable.
Hence the assertions $(i) $ and $(ii)$ hold. \hfill$\Box$
\begin{Prop}
Let  be the  fractional Lagrange system $(3.3) $ and $~q\in (0,1).$\\
$(i) $ If $~c_{1} < 0, c_{2} < 0,~$ then $~ e_{2}^{m}~$ and $~ e_{3}^{m}~$ are asymptotically stable for all\\
$~m\in ( -\sqrt {c_{1} c_{2}},  \sqrt {c_{1} c_{2}}) ~$ and unstable  $~(\forall) m\in (-\infty,  -\sqrt {c_{1} c_{2}})\cup  (\sqrt {c_{1} c_{2}}, \infty).$ \\
$(ii)~$  If $~c_{1} > 0,~c_{2} < 0~$  or $~ c_{2} > 0,~ c_{1}\in {\bf R} ,$ then $~ e_{2}^{m}~$ and $~ e_{3}^{m}~$ are unstable  $~(\forall) m \in  {\bf R}^{\ast}.~$
\end{Prop}
{\it Proof.} We consider $~k\in \{2, 3\}.~$ The Jacobian matrices of  $ (3.3)$ at $~e_{k}^{m}~$ are\\
$~J(e_{2}^{m}, c_{1}, c_{2})=\left (
\begin{array}{ccccc}
c_{1} & 0 & m \\[0.2cm]
0 & c_{2} & 0 \\
m  & 0 &  c_{2}\\
\end{array}\right )~ $ and $~J(e_{3}^{m}, c_{1}, c_{2})=\left (
\begin{array}{ccccc}
c_{1} & m& 0 \\[0.2cm]
m & c_{2} & 0 \\
0 & 0 &  c_{2}\\
\end{array}\right ).$
Its characteristic polynomials are $~p_{k}(\lambda) = -(\lambda -
c_{2})[ \lambda^{2} - (c_{1} + c_{2})\lambda + c_{1}c_{2}- m^{2}].~$  The roots of the equation $~ p_{k}(\lambda)= 0~$ are $ \lambda_{1}=
c_{2},~\lambda_{2,3}= \frac{(c_{1}+c_{2})\pm \sqrt{\Delta}}{2}\in {\bf R},$ where $~\Delta = (c_{1}- c_{2})^{2} +4 m^{2}.$\\
$(i)~$ We suppose  $~c_{1} <0, c_{2} < 0. $ We have $~\lambda_{2}, \lambda_{3}< 0~$ if and only if $~\lambda_{2} + \lambda_{3}< 0~$ and $~\lambda_{2}\lambda_{3}> 0~$. Then $~c_{1}+c_{2} < 0~ $ and $~(c_{1} + c_{2})^{2}- \Delta >0.~$
 It follows that $~\lambda_{i} < 0 , i=\overline{1,3}~$ for all $~ m\in (-\sqrt {c_{1} c_{2}},  \sqrt {c_{1} c_{2}}).$  Therefore,  $~ e_{k}^{m}~$ is  asymptotically stable.\\
$ (ii)~$ We suppose  $~c_{1} > 0~$ and $~c_{2} < 0~$  or $~ c_{2} > 0~$ and $ c_{1}\in {\bf R}.$ In this cases,$ ~ J(e_{k}^{m}, c_{1}, c_{2})~$ has at least a positive eigenvalue and so $~e_{k}^{m}~$ is unstable $~(\forall) m\in {\bf R}^{\ast}~.$ Therefore, the assertions $ (i) $ and $ (ii) $ hold.\hfill$\Box$
\begin{Ex}
{\rm By choosing the control parameters $ c_{1}, c_{2} $ that satisfy one condition from Proposition 3.3, then the
trajectories of the fractional-order Lagrange system $(3.3)$ are driven to the stable equilibrium point $e_{3}^{m}~(m\neq 0)$.
 For example, we select $ c_{1} =-1.75, c_{2} = -2,~$ then the stability condition $(i)$ of Proposition 3.3 is achieved. This
implies that, the trajectories of the system $(3.3)$ converge to the equilibrium point $e_{3}= ( 0, 0, m)$ when $ m\in (-1.8708, 1.8708)~$ and $~q\in (0, 1).$ For example, the  fractional-order Lagrange system is asymptotically stable at $ e_{3}=(0, 0, 1.75)$
for $~q \in (0,1)$.}\hfill$\Box$
\end{Ex}
Using Matlab, in Table 3.1 we give a set of values for $~c_{1}, c_{2},~$  the equilibrium states and corresponding eigenvalues of fractional-order Lagrange system $(3.3).$
{\small
$$\begin{array}{|l|l|l|l|l|} \hline
       c_{1}, c_{2}\in {\bf R}^{\ast}   & Eigenvalues & m,~q &  e_{i}^{m}  & Stability  \\ \hline

 c_{1}=-0.2,~c_{2}=-0.8  & \lambda_{1}=-0.2,~\lambda_{2,3}=-0.8 &  q\in (0,1) & e_{0} &  stable \cr \hline

 c_{1}=-2, ~ c_{2}=-1.85 & \lambda_{1}=-2,~\lambda_{2}=-0.85 & m=1 & e_{1}^{m} &  stable \cr
  & \lambda_{3}=-2.85 & q\in (0,1) &  & m\in (-1.85, 1.85) \cr \hline

 c_{1}=-7.2, ~ c_{2}=-0.2 & \lambda_{1}=-0.2,~\lambda_{2}=0.4182 & m=-2 & e_{2}^{m} &  unstable \cr
  & \lambda_{3}=-7.8182 & q\in (0,1) &  & m\in (-\infty, -1.2)\cup (1.2, \infty) \cr \hline

 c_{1}=-1.75, ~ c_{2}=-2 & \lambda_{1}=-2,~\lambda_{2}=-0.9911 & m=1.75 & e_{3}^{m} &  stable \cr
  & \lambda_{3}=-2.7588 & q\in (0,1) &  & m\in (-1.8708, 1.8708) \cr \hline
\end{array} $$\\[-0.4 cm]
{\bf Table 3.1.} {\it The controls $~c_{1}, c_{2},~$ equilibrium states $~e_{i}^{m}~$ and corresponding eigenvalues}.}\\[-0.5cm]

\section{Numerical integration of the fractional-order Lagrange system $(3.3)$}

Consider the following general form of the initial value problem (IVP) with Caputo derivative \cite{odmo}:\\[-0.4cm]
\begin{equation}
D_{t}^{q} y(t) = f(t,y(t)),~~~ y(0)=y_{0},~~~t\in I=[0,T],~T>0  \label{(4.1)}
\end{equation}
where $~y: I \rightarrow {\bf R}^{n},~f: {\bf R}^{n} \rightarrow {\bf R}^{n}~$ is a continuous nonlinear function and
$ q\in (0,1), $ represents the order of the derivative.

The right-hand side of the IVP $ (4.1) $ in considered examples are Lipschitz functions and the numerical method used in this works to integrate system $ (4.1) $ is the {\it Fractional Euler's method}.

Since $ f $ is assumed to be continuous function, every solution of the initial value problem given by $ (4.1) $ is also a solution of the following {\it Volterra fractional integral equation}:\\[-0.4cm]
\begin{equation}
y(t) = y(0) +~ I_{t}^{q}f(t,y(t)), \label{(4.2)}
\end{equation}
where $ I_{t}^{q} $ is the $ q-$order Riemann-Liouville integral operator $ (2.1). $\\[-0.4cm]
\begin{equation}
I_{t}^{q}f(t) =\displaystyle\frac{1}{\Gamma(q)}\int_{0}^{t}{(t-s)^{q
-1}}f(s,y(s))ds,~q > 0,~~~s\in [0, T].  \label{(4.3)}
\end{equation}
Moreover, every solution of $ (4.2) $ is a solution of the (IVP) $ (4.1).$

To integrate the fractional equation $( 4.1),$ means to find the solution of $ (4.2) $ over the interval $~[0,T]. $ In this context, a set of points $~(t_{j}, y(t_{j})) $ are produced which are used as approximated values. In order to achieve this approximation, the interval $ [0,T] $  is partitioned into $ n $ subintervals $~[t_{j}, t_{j+1} ] $  each equal width $~ h =\frac{T}{n}, ~ t_{j} = j h $  for  $~ j = 0,1,..., n.~ $ For the fractional-order $ q $ and $~ j = 0, 1, 2, ... ,$ it computes an approximation denoted as $~y_{j+1} ~$ for $~ y ( t_{j+1}),~ j= 0, 1, ... .$

The general formula of the fractional Euler's method for to compute the elements $~y_{j},$ is\\[-0.4cm]
\begin{equation}
y_{j+1} = y_{j} + \displaystyle\frac{h^{q}}{\Gamma(q+1)} f( t_{j}, y(t_{j})),~~~~~ t_{j+1} = t_{j} + h,~~~ j=0, 1, ..., n. \label{(4.4)}
\end{equation}
For more details, see \cite{odmo, danc, ahme}.

For the numerical integration of the system $ (3.3),$ we apply the fractional Euler method (FEM). For this, consider the following fractional differential equations\\[-0.2cm]
\begin{equation}
\left\{\begin{array}{lcl}
 D_{t}^{q} x^{i}(t) & = &
F_{i}(x^{1}(t), x^{2}(t), x^{3}(t)),~~~ i=\overline{1,3}, ~~t\in
(t_{0}, \tau),~q \in (0,1)\\
x(t_{0}) &=& (x^{1}(t_{0}), x^{2}(t_{0}), x^{3}(t_{0}))
\end{array}\right.\label{(4.5)}\\[-0.2cm]
\end{equation}
where\\[-0.4cm]
\begin{equation}
\left\{\begin{array}{lcl} F_{1}(x(t)) & = &x^{2}(t) x^{3}(t) + c_{1} x^{1}(t),\\
F_{2}(x(t)) & = & x^{1}(t) x^{3}(t) + c_{2}x^{2}(t),~~~~~  c_{1}, c_{2}\in {\bf R}^{*},\\
F_{3}(x(t)) & = & x^{1}(t) x^{2}(t) + c_{2} x^{3}(t).\\
\end{array}\right.\label{(4.6)}
\end{equation}
 Since the functions $ F_{i}(x(t)), i=\overline{1,3} $ are continuous, the initial value problem $(4.5)$ is equivalent
 to system of Volterra integral equations, which is given as follows:\\[-0.2cm]
\begin{equation}
x^{i}(t)~=~ x^{i}(0)  + ~I_{t}^{q} F_{i}(x^{1}(t), x^{2}(t), x^{3}(t)),
~~~~~i=\overline{1,3}.\label{(4.7)}
\end{equation}
The system $(4.7)$ is called the {\it Volterra integral equations
associated to fractional-order Lagrange system} $(4.5)$.

The problem for solving the system $(4.5)$ is reduced to one of solving a sequence of systems of fractional equations in
increasing dimension on successive intervals $[j, (j+1)]$.

For the numerical integration of the system $(4.6)$ one can use the fractional Euler method (the formula $(4.4)$ ), which is expressed as follows:\\[-0.4cm]
\begin{equation}
x^{i}(j+1)=x^{i}(j)+ \frac{h^{q}}{\Gamma (q+1)} F_{i}(x^{1}(j), x^{2}(j),
x^{3}(j)),~~~ i=\overline{1,3},\label{(4.8)}
\end{equation}
where $ j=0,1,2,...,N,  h=\displaystyle\frac{T}{N}, T>0, N>0.~$

More precisely, the numerical integration of the fractional system $(4.5)$ is given by:\\[-0.4cm]
\begin{equation}
\left \{ \begin{array}{ll} x^{1}(j+1) &= x^{1}(j)+
h^{q}~\displaystyle\frac{1}{\Gamma(q+1)}(
 x^{2}(j) x^{3}(j) + c_{1} x^{1}(j))\\[0.3cm]
x^{2}(j+1) &= x^{2}(j)+
h^{q}~\displaystyle\frac{1}{\Gamma(q+1)}( x^{1}(j) x^{3}(j) + c_{2} x^{2}(j))\\[0.3cm]
 x^{3}(j+1) &= x^{3}(j) +
h^{q}~\displaystyle\frac{1}{\Gamma(q+1)}(
 x^{1}(j) x^{2}(j) + c_{2}x^{3}(j))\\[0.3cm]
x^{i}(0)&= x_{e}^{i}+\varepsilon,~~~i=\overline{1,3}.\\
\end{array}\right. \label{(4.8)}
\end{equation}
Using \cite{odmo, diet}, we have that the numerical algorithm given by
$(4.9)$ is convergent.\\[-0.2cm]
\begin{Ex}
{\rm  Let us we present the numerical simulation of solutions of the fractional-order Lagrange system  $(3.3)$ which has considered in Example $~3.1.$\\
For this we apply the algorithm $ (4.9) $ and software Maple. Then, in $(4.9)$ we take: $~ c_{1}= -1.75, c_{2}= -2,~h = 0.01, \varepsilon= 0.01,  N = 500, t = 502. $

The orbits $(n, x^{i}(n)), i=\overline{1,3}~$  for the solutions of fractional-order Lagrange system for the equilibrium state $ e_{3} = (0, 0, 1.75)~ $ have the representations given in figures Fig. 1(a), 2(a), 3(a) (for $~q = 0.65 ~)$ and Fig. 1(b), 2(b), 3(b) (for $~ q = 1).~$}\hfill$\Box$
\end{Ex}
\begin{center}
\begin{tabular}{ccc}
\includegraphics[width=5cm]{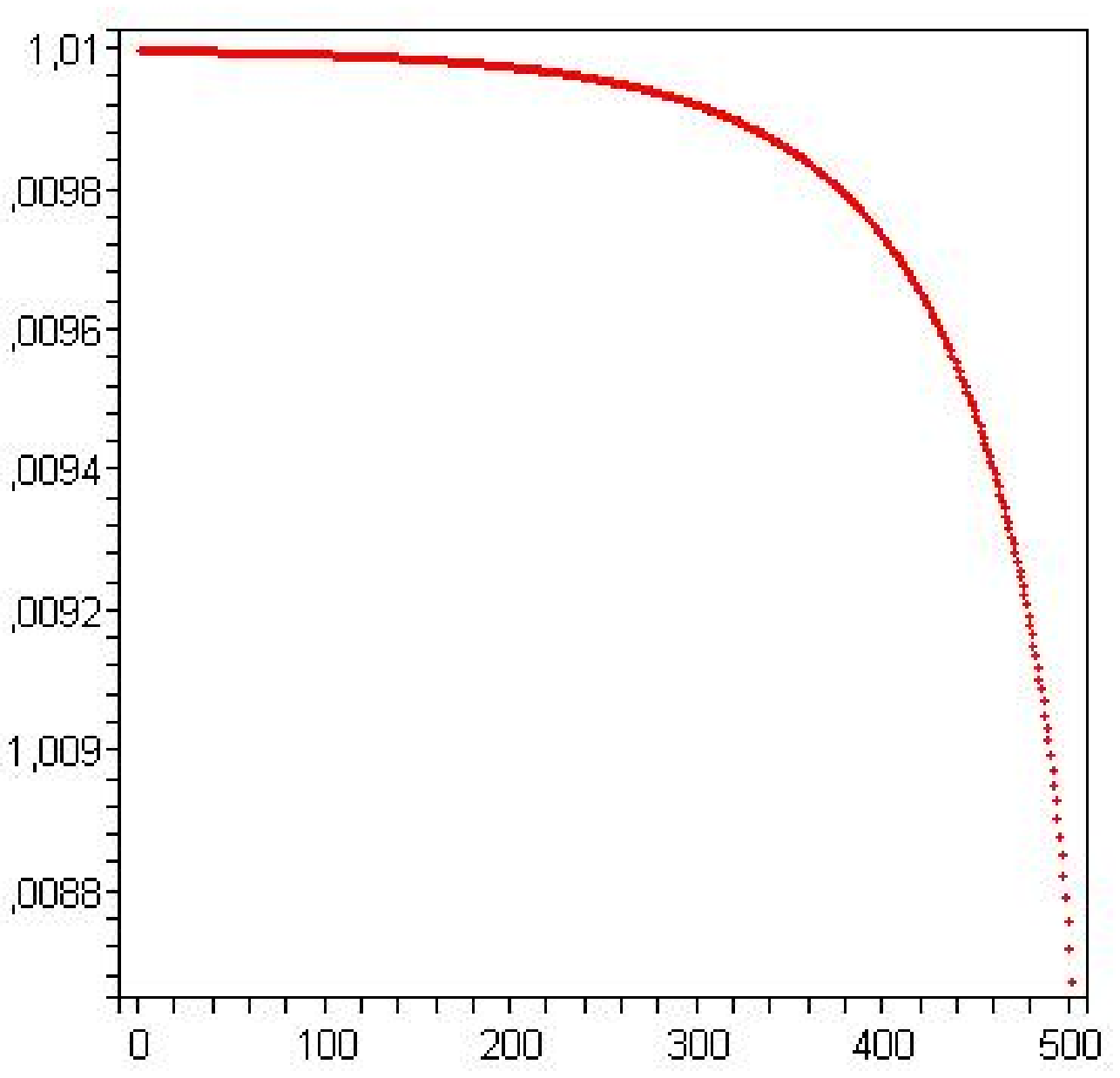}&&
\includegraphics[width=5cm]{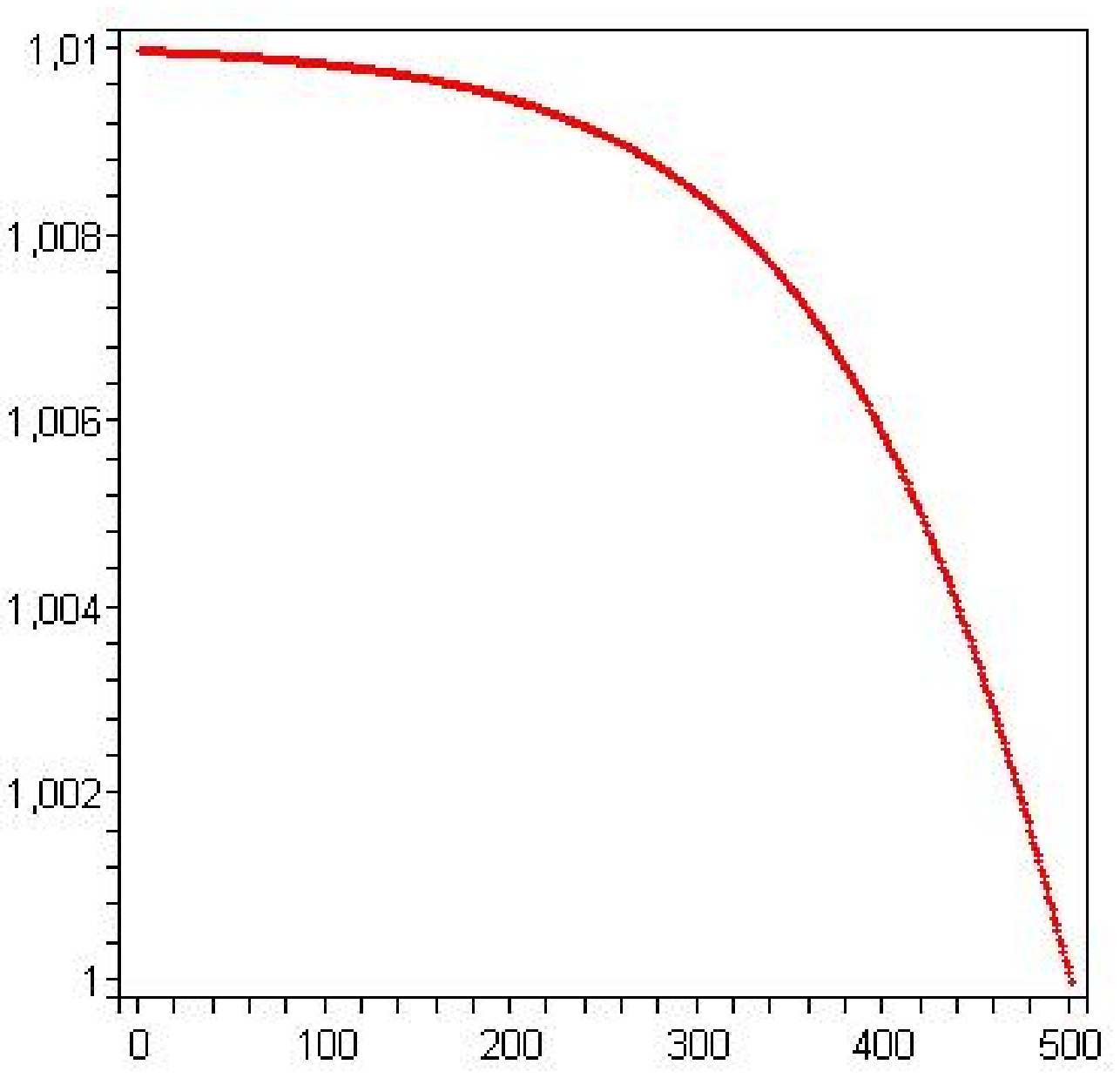}
\end{tabular}
\end{center}
\begin{center}
\begin{tabular}{ccc}
{\bf Fig.1.(a)} $~(n, x^{1}(n))$ for $ q = 0.65 $ & & {\bf Fig.1.(b)}
$~(n, x^{1}(n))$ for $ q =1 $
\end{tabular}
\end{center}
\begin{center}
\begin{tabular}{ccc}
\includegraphics[width=5cm]{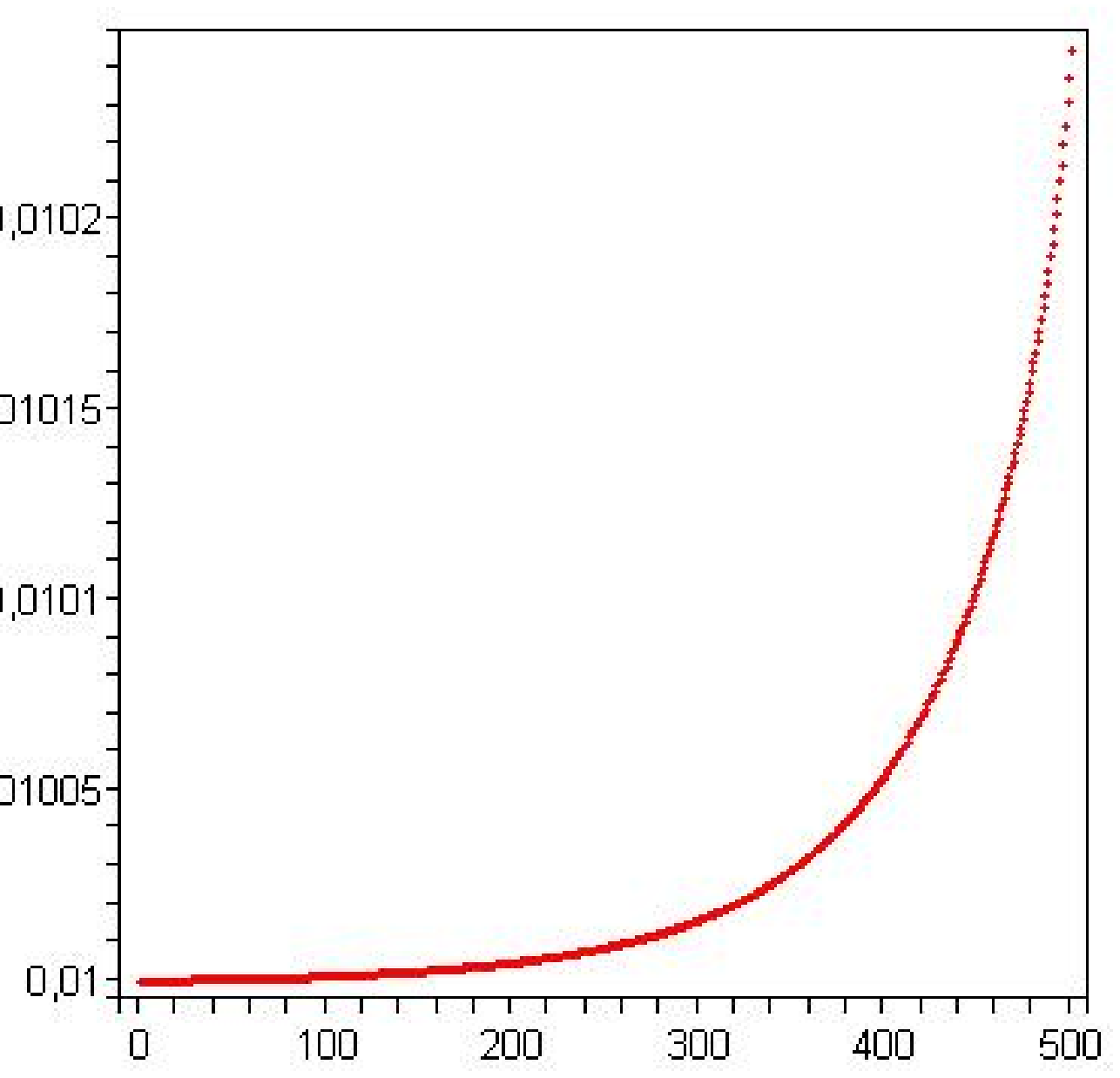}& &
\includegraphics[width=5cm]{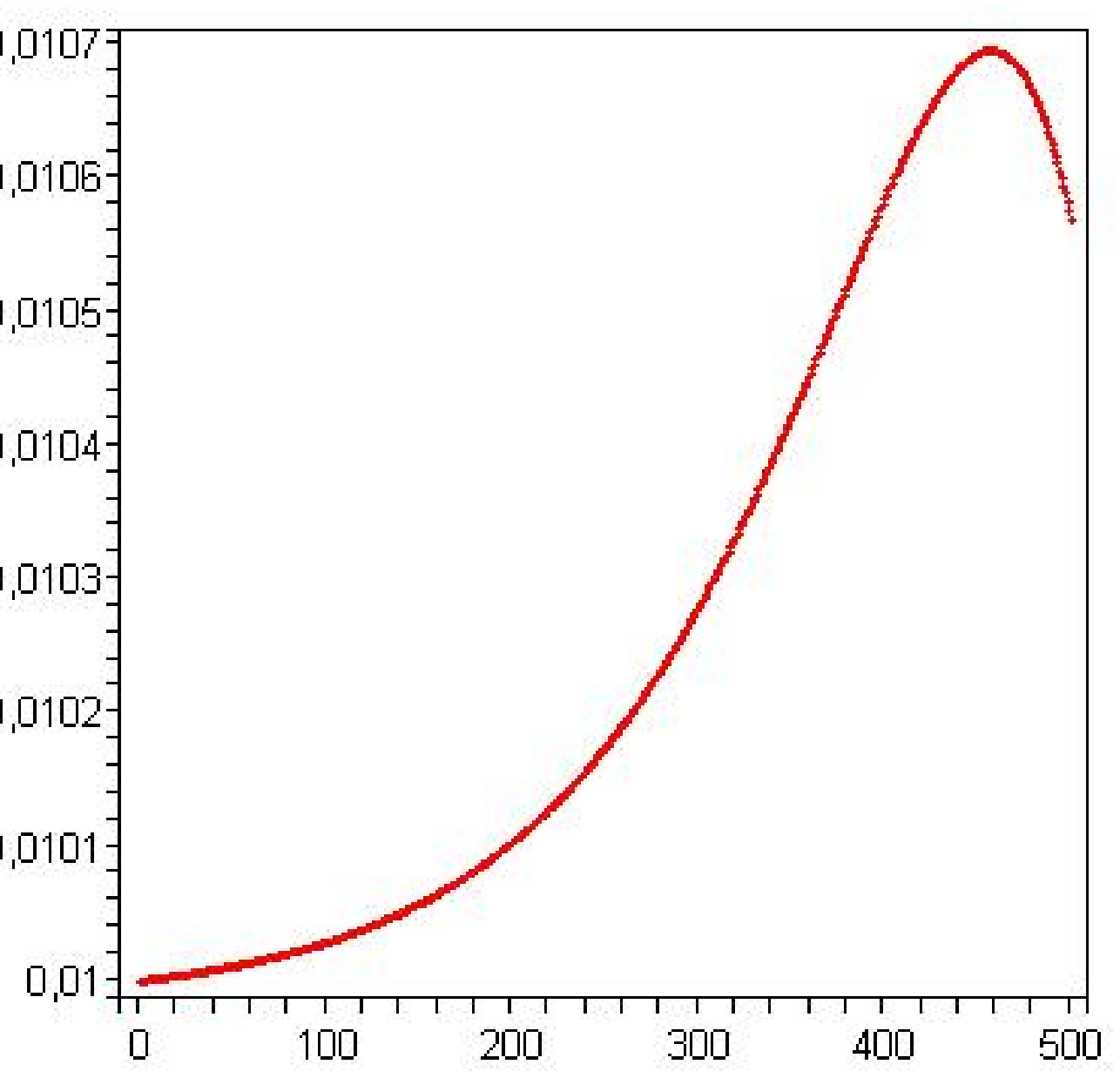}
\end{tabular}
\end{center}
\begin{center}
\begin{tabular}{ccc}
{\bf Fig.2.(a)} $~(n, x^{2}(n))$ for $ q = 0.65 $ & & {\bf Fig.2.(b)}
$~(n, x^{2}(n))$ for $ q = 1 $
\end{tabular}
\end{center}
\begin{center}
\begin{tabular}{ccc}
\includegraphics[width=5cm]{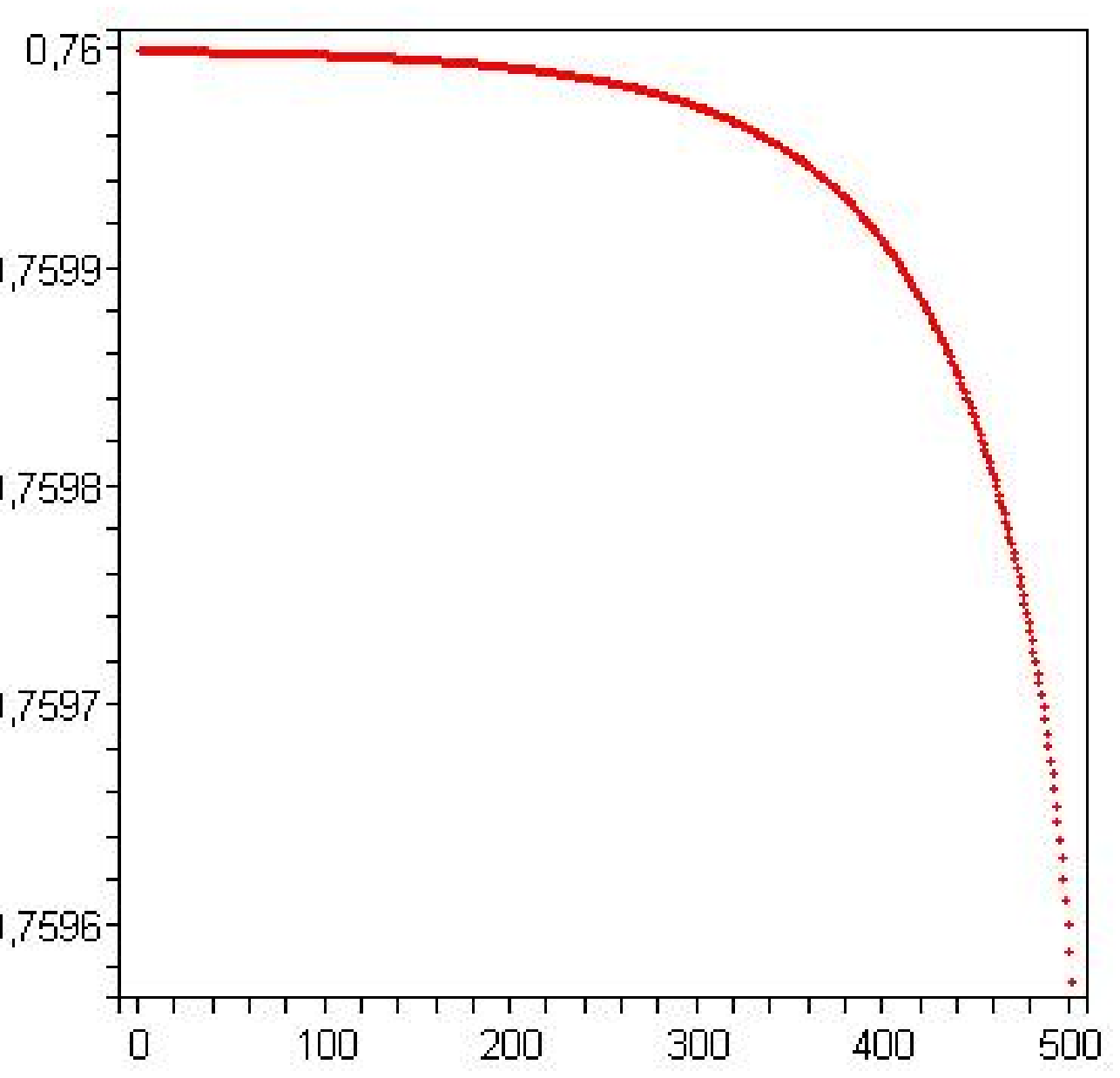}& &
\includegraphics[width=5cm]{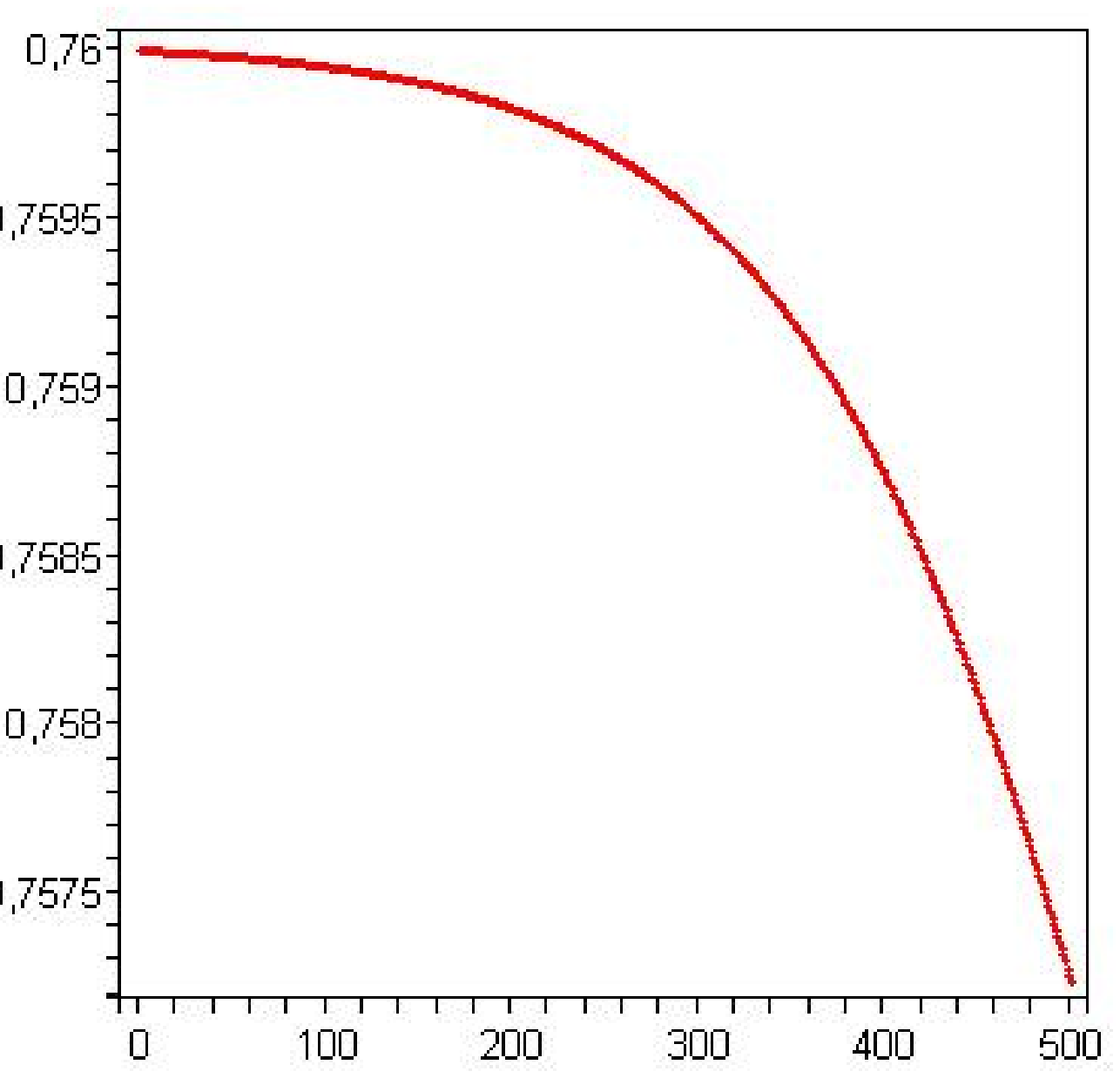}
\end{tabular}
\end{center}
\begin{center}
\begin{tabular}{ccc}
{\bf Fig.3.(a)} $~(n, x^{3}(n))$ for $ q = 0.65 $ & & {\bf Fig.3.(b)}
$~(n, x^{3}(n))$ for $ q = 1 $
\end{tabular}
\end{center}
The numerical simulations confirm the validity of the theoretical analysis.\\[-0.4cm]
\begin{Rem}
{\rm  Appyling  $ (4. 9) $ and  Maple for the numerical simulation of solutions of fractional Lagrange system $~(3.3)~$ for each pair $~(c_{1}, c_{2})~$ of values given in the Table (3.1), it will be found that the results obtained are valid.}~~~\hfill$\Box$
\end{Rem}
{\bf Conclusions.} This paper presents the fractional-order Lagrange system  $(3.3) $ associated to  system $(1.1).$  The fractional-order Lagrange system $ (3.3) $ was studied from fractional differential equations theory point of view: asymptotic stability, determining of sufficient conditions on parameters $ c_{1}, c_{2} $ to control the chaos in the proposed fractional system and numerical integration of the fractional model $ (4.2).~$  By choosing the right parameters $ c_{1} $ and $ c_{2}$ in the fractional model  $~(3.3),~$ this work offers a series of chaotic fractional differential systems. \hfill$\Box$

{\bf Acknowledgments.} The author has very grateful to be reviewers for their comments and suggestions. \\[-0.5cm]

Author's adress\\[-0.5cm]

Mihai Ivan\\[0.1cm]
West University of Timi\c soara. Seminarul de Geometrie \c si Topologie.\\
Teacher Training Department. Timi\c soara, Romania.\\
E-mail: mihai.ivan@e-uvt.ro\\

\begin{thebibliography}{99}
%\addcontentsline{toc}{chapter}{Bibliografie}
\smallskip
\bibitem{agra} O.P. Agrawal, {\it Solution for a fractional diffusion-wave equation defined in a bounded domain}, Nonlinear Dynamics, {\bf 29} (2002), no. 1–4, 145–155.

\bibitem{ahme} H.M. Ahmed {\it Fractional Euler method: An effective tool for solving fractional differential equations of fractional order}, J. of Egyptian Mathematical Society, {\bf 26} (2018), no.1, 38-43.

\bibitem{ahma} W.M. Ahmad, J.C. Sprott, {\it Chaos in fractional order autonomous nonlinear systems}, Chaos Solitons Fractals {\bf 16}(2003), 339–351.

\bibitem{badi} D. B\u aleanu, K. Diethelm, E. Scalas, J.J. Trujillo, {\it Fractional Calculus Models and
Numerical Methods}, Series on Complexity, Nonlinearity and Chaos. World Scientific, 2012.

\bibitem{bhad} S. Bhalekar, V. Daftardar-Gejji, {\it Synchronization of different fractional order chaotic systems using active control}, Communications in Nonlinear Science and Numerical Simulation, {\bf 15} (2010), no. 11,  3536–3546.

\bibitem{danc} M.F. Danca, {\it Hidden chaotic attractors in fractional-order systems}, Nonlinear Dynamics, {\bf 89} (2017), no. 1, 577-586.

\bibitem{demi}  M. Degeratu, M. Ivan, {\it Linear connections on Lie algebroids}, Proceed. of the 5th Conf. of Balkan Society of Geometers, Aug. 29- Sept. 2, 2005. Mangalia, Romania. Geometry Balkan Press, 2006, 44-53.

\bibitem{difo} K. Diethelm, N.J. Ford,  {\it Analysis of fractional differential equations}, J. Math. Analysis Applications, {\bf 265}(2002), 229–248.

\bibitem{diet} K. Diethelm, {\it  The Analysis of Fractional Differential Equations: An Application-Oriented Exposition  Using Differential Operators of Caputo Type}, Springer, 2010.

\bibitem{nabu} R.A. El-Nabulsi, {\it A fractional action-like variational approach of some classical quantum and geometrical dynamics}, Int. J. Appl. Math.,{\bf 17}(2005), 299-317.

\bibitem{enpp} R.D. Ene, C. Pop,  C. Petri\c sor, {\it Systematic review of geometric approaches and analytical integration for Chen’s system}, Mathematics, {\bf 8} (2020),no. 9,  1530.

\bibitem{gunu} H. G$\ddot{u}$mral, Y. Nutku, {\it Poisson structure of dynamical systems with three degrees of freedom}, J. Math. Phys. {\bf 34} (1993), no. 12.

\bibitem{giv1} G. Ivan, {\it Geometrical and dynamical properties of general Euler top system}, Indian J. Pure Appl. Math., {\bf 44} (2013), no.1, 77-93.

\bibitem{giv2} G. Ivan, {\it On fractional differential equations of 3D Maxwell-Bloch type}, Int. J. Geom. Method in Modern Physics, {\bf 11}(2014), no. 4, 1450028 (12 pages).  Doi:10.1142/S0219887814500285.

\bibitem{giop} G. Ivan, D. Opri\c s, {\it Dynamical systems on Leibniz algebroids}, Diff. Geometry-Dynamical systems, {\bf 8}(2006). Geometry Balkan Press, 127-137.

\bibitem{gimo}  G. Ivan, M. Ivan, D. Opri\c s, {\it Fractional Euler-Lagrange and fractional Wong equations for Lie algebroids}. Proceed. of The 4-th Int. Colloq. "Math. in Eng. and Numerical Phys." October 6-8, 2006, Bucharest, Romania, 73-80.[Google Scholar]
\bibitem{imod} G. Ivan, M. Ivan, D. Opri\c s, {\it Fractional dynamical systems on fractional Leibniz algebroids}. Analele \c Stiin\c t.Univ, "Al. I. Cuza" din Ia\c si (S.N.), Matematic\u a, {\bf 53} (2007), Supl., 222-234.

\bibitem{igmp} G. Ivan, M. Ivan, C. Pop, {\it Numerical integration and synchronization for the 3-dimensional metriplectic Volterra system}, Mathematical Problems in Engineering, {\bf 2011}, Article ID 723629 (11 pages). Doi:10.1155/2011/723629.

\bibitem{opmi} G. Ivan, D. Opri\c s,  M. Ivan, {\it Stochastic fractional equations associated to Euler top system}. Int. J. Geom. Met. Mod. Phys., {\bf 10}(2013), no. 1, 1220018 (10 pages).\\  Doi:10.1142/S02198878.12200186.

\bibitem{ivmi} M. Ivan, {\it Control chaos in the fractional Lorenz-Hamilton system}, Fractional Differential Calculus, {\bf 6}(2016), no. 1, 111-119. Doi:10.7153/fdc-06-07.

\bibitem{migi} M. Ivan, G. Ivan, {\it On the fractional Euler top system with two parameters}.  Int. J. Modern Eng. Research, {\bf 8}(2018), no. 4, 10-22

\bibitem{migo} M. Ivan, G. Ivan, D. Opri\c s, {\it Fractional equations of the rigid body on the pseudo-orthogonal group $ SO(2,1) $}, Int. J. Geom. Method in Modern Physics, {\bf 6}(2009), no. 7, 1181-1192.

\bibitem{laph} C. L\u azureanu, C. Petri\c sor, C. Hedrea, {\it On a deformed version of the two-disk dynamo}, Appl. of Math. {\bf 66}(2021),no. 3, 345-372.

\bibitem{lima} J. Llibre, Y.P. Martinez, {\it Dynamics of a family of Lotka-Volterra Systems in $~{\bf R}^{\ast}, $} Nonlinear Anal.{\bf 199}(2020),111915.

\bibitem{mati} D. Matignon, {\it Stability results for fractional differential equations with applications to control processing}. In: Proceedings of the Computational Engineering in Systems and Applications, IMACS, IEEE-SMC, Lille, France, July 1996, {\bf 2}(1996), 963-968.

\bibitem{odmo} Z.M. Odibat, S. Momani, {\it An algorithm for the numerical solution of differential equations of fractional order}, J. Appl. Math. \& Informatics, {\bf 26} (2008), 15-27.

\bibitem{podl} I. Podlubny, {\it Fractional Differential Equations}, Academic Press, 1999.

\bibitem{pagi} C. Pop, A. Aron, C. Galea, M. Ciobanu, M. Ivan, {\it Some geometric aspects in theory of Lotka-Volterra system}, Proc. of the 11-th WSEAS
Int. Conf. on Sust. in Science Engineering, Timi\c soara, Romania, May (2009), 91-97.[Google Scholar]

\bibitem{puta} M. Puta, C. Pop, C. Danaiasa, C. Hedrea, {\it Some geometric aspects in the theory of Lagrange system}, Tensor, N.S. {\bf 69}(2008), 83-87.

\bibitem{takh} L. Takhtajan, {\it On foundation of the generalized Nambu mechanics}, Communications in  Mathematical Physics, {\bf  160} (1994), no. 2, 295–315.\\[-0.5cm]
\end{thebibliography}
\end{document}